\documentclass [10pt]{article}\input amssym.def \input amssym

\begin {document}
\def\cir{\circle{5}}
\def \Z{\bf Z}
\def \N{\bf N}
\newtheorem{theorem}{\bf Theorem}
\newtheorem{definition}{\bf Definition}
\newtheorem{example}{\bf Example}
\newtheorem{remark}{\bf Remark}

 \author {Eduardo Saenz de Cabezon
 \\ Universidad de La Rioja
 \\ eduardo.saenz-de-cabezon@dmc.unirioja.es} 
 \title {Bordered magic squares: Elements for a comprehensive approach \footnote {Subject Mathematics Classification 00A08,97A20,05B15}} 
 \date {Oct-2003} \maketitle

\begin {abstract}
General methods for the construction of magic squares of any order
have been searched for centuries. Several ``standard
strategies'' have been found for this purpose, such as the ``knight movement'', or the
construction of bordered magic squares, which played an important role
in the development of general methods. 

What we try to do here is to give a general and comprehensive
approach to the construction of magic borders, capable of assuming
methods produced in the past like particular cases. This
general approach consists of a transformation of the problem of
constructing magic borders to a simpler - almost trivial - form.
In the first section, we give some definitions and notation. The second section consists of
the exposition and proof of our method for the different
cases that appear (Theorems 1 and 2). As an application of this method, in the third section we caracterise magic borders of even order, giving therefore a first general result  for bordered magic squares.
 
 Although methods for the construction of bordered magic squares have always 
 been presented as individual succesful attempts to solve the 
 problem, we will see that a common pattern underlies the 
 fundamental mechanisms that lead to the construction of such 
 squares. 
 This approach provides techniques for constructing many magic 
 bordered squares of any order, which is a first step to construct all 
 of them,and finally know how many bordered squares are for any 
 order. These may be the first elements of a general theory on 
 bordered magic squares.
\end {abstract}

\pagebreak
 
\section {Previous considerations}
Let us begin with some useful definitions:

\begin {definition}
{\it Magic square of order $n$:} An $n\times n$ array 
containing the integers $1$ through $n^2$, each one once, arranged in 
such a way that all rows, columns and the two main diagonals have the 
same sum $S_{n}=\frac {n}{2} (n^2+1)$, which is called its magic 
constant.
\end{definition}

\begin {definition} Border of a magic square: The set consisting of its first and 
last row and column.\end{definition}

\begin {definition} Complementaries: In a magic square of order 
$n$, a pair of numbers $a,b$ so that $a+b=n^2+1$ we will denote 
$\bar \alpha$ the complementary of $\alpha$\end{definition}

\begin {definition} Magic border: A border that verifies that 
every line (row or column) sums $S_{n}$ and that every two opposite 
numbers are complementaries. Thus, a magic border consists of tho 
subsets, one containing the complementaries of the other, which will 
be denoted $A$ and $A'$ respectively.\end{definition}

\section { A method for the construction of magic borders}

Since the beginning of the history of bordered magic squares, many particular methods for their construction have been produced, references on this topic can be seen at \cite{Ale}\cite{Descombes}\cite{Fre}\cite{Ses}\cite{Suz}\cite{Tra}\cite{Wei}. What we try to do here is to give a general approach to these methods.

\subsection{Exposition of the general method for the construction of Magic Borders} 
As we said before, the method we propose consists
of a simplification of the problem of building magic borders. We
will put the problem in terms that will lead us to an easy way of
resolution.

First of all, we name the numbers in the border as follows:
$$\matrix{v&b_{1}&b_{2}&\dots&b_{n}&w \cr
c_{1}&\quad&\quad&\quad&\quad&\bar c_{1}\cr
c_{2}&\quad&\quad&\quad&\quad&\bar c_{2}\cr
\dots&\quad&\quad&\quad&\quad&\dots\cr
c_{n}&\quad&\quad&\quad&\quad&\bar c_{n}\cr
\bar w&\bar b_{1}&\bar b_{2}&\dots&\bar b_{n}&\bar v}$$

\noindent where $\bar \alpha$ is the {\it complementary} of $\alpha$. Note that
$\alpha + \bar \alpha=(n+2)^2+1=n^2+4n+5$, which is the difference between the
magic constant of the big square and that of the inner one (The magic constant
 of a magic square of order $n$ is
$S_{n}=\frac{n}{2}(n^2+1)$. So, the constant of the big square is
$S_{n+2}=\frac{(n+2)}{2}((n+2)^2+1)$ and the constant of the inner
square  ---with it's numbers incremented in $2n+2$---  is
$S_{n}+n(2n+2)$. Hence, $S_{n+2}-S_{n}-n(2n+2)=n^2+4n+5$). So, we
see that every row, column, and diagonal sums the same, the
magic constant, except the first and last row and column. We
just need then to verify the magic condition placing the numbers
we have: $\{ 1,\ldots ,2n+2\} \cup \{ n^2+2n+3, \ldots ,(n+2)^2 \}$
so that the following equalities stand
 $$ \sum_{i=1}^nb_{i}+v+w=\sum_{i=1}^n
\bar b_{i}+\bar v+\bar w=\sum_{i=1}^nb_{c}+v+\bar w=\sum_{i=1}^n c_{i}+\bar v+w=S_{n+2}$$
 To
make the placement of these numbers easier, we arrange them in two
columns, as follows

\bigskip

\begin {picture}(60,60)(-80,0) \put(0,60){$1$}
\put(10,60){$\bullet$}
 \put(40,60){$\bullet$} \put(50,60){$(n+2)^2$}
 \put(0,45){$2$} \put(10,45){$\bullet$}
 \put(40,45){$\bullet$} \put(50,45){$(n+2)^2-1$}
 \put(12,30){$\vdots$}\put(42,30){$\vdots$}
\put(12,15){$\vdots$}\put(42,15){$\vdots$}
\put(0,0){\makebox(0,0)[br]{$2n+2$}}\put(10,0){$\bullet$}
 \put(40,0){$\bullet$} \put(50,0){$n^2+2n+3$}
 \end{picture}
\medskip

 With this arrangement, if we match two numbers with a line, say
 $x$ and $y$, we can assign to that matching line, the number that
 expresses the difference between the sum of those two numbers and
  the sum of two complementary numbers, i.e. $x+y-(n^2+4n+5)$. So,
  for example

  \begin {picture}(160,60)(-80,15) \put(0,50){$x$}
\put(10,50){$\bullet$}
 \put(40,50){$\bullet$}
  \put(10,40){$\bullet$}
 \put(40,40){$\bullet$}
 \put(10,30){$\bullet$}
 \put(40,30){$\bullet$} \put(50,30){$y$}
 \put(12,52){\line(3,-2){30}} \put (60,40){equals $-2$}
 \put(150,50){$\bullet$}
 \put(180,50){$\bullet$}\put(190,50){$y$}
  \put(150,40){$\bullet$}
 \put(180,40){$\bullet$}
 \put(150,30){$\bullet$}
 \put(180,30){$\bullet$} \put(140,30){$x$}
 \put(152,32){\line(3,2){30}} \put(200,40){equals $+2$}
 \end{picture}

 We denote this number by $d_{xy}$.

 Next, we will show different procedures for the construction of the
 magic border if the square has even or odd order.

 --- Even order:

 If the order of the square is even, we proceed as follows:

 We choose two numbers $v$ and $\bar w$ to be placed in the corners \footnote {In section 3 we will discuss what numbers can or cannot be chosen to be $v$ and $w$}. After this,
 we take $n$ numbers from one column and the other alternately,
 beginning in the opposite column to that containing the number
 named $v$, and taking care not to have chosen any pair of
 complementary numbers. These will be our $b_i 's$. The way of
 choosing these $b_i 's$ has to be one in which we have some 
 permutation $\sigma$ of $\{1,\dots,n\}$ and a way to put
 them in pairs, $$\beta = \{ (b_{\sigma (1)},v),(b_{\sigma (2)},b_{\sigma (3)}),
 \dots,(b_{\sigma (n)},w)\}$$ verifying $$\sum_{(x,y)\in\beta}
 d_{xy}=0.$$
After this, we have to choose $n$ numbers, from one column and the
other alternately, in such a way that not two of them are 
complementary.
These will be the $c_i 's$ and they have to be so that there
exists some permutation $\tau$ of $\{1,\dots,n\}$ and a way
 to put them in pairs, $$\gamma = \{ (c_{\tau(1)},c_{\tau(2)}),
 \dots,(c_{\tau(n-1)},c_{\tau(n)})\}$$ verifying
$$\sum_{(x,y)\in\gamma}
 d_{xy}=-d_{v\bar w}.$$

 Making such an election of the numbers,we assure that putting
 $b_i's$ in the top (or bottom) row, except the corners, in which we will put
 $v$ and $w$; and putting $c_i 's$ in the column between $v$ and
 $\bar w$ the border obtained is magic, regardless the order in which
 $b_i's$ and $c_i 's$ have been placed in their respective rows or
 columns.

 --- Odd order:

 If the order of the square is odd, the procedure is a little different.
 First, we choose a number for one of the squares, $v$, and
 denote $d_{v}:=v-\frac{(n+2)^2+1}{2}$. After that, as we have done in the
 even case, we choose $n+1$ numbers (which will be our $b_{i}'s$ and
 $w$)
 from one column and the other, alternately, in such a way that we can
 match them in pairs, being $\beta = \{ (b_{\sigma (1)},v),(b_{\sigma (2)},
 b_{\sigma (3)}),\dots,(b_{\sigma (n)},w)\}$ 
 verifying $$\sum_{(x,y)\in\beta}d_{xy}+d_{v}=0.$$
 After this, we have to choose $n$ numbers, from one column and the
other alternately, in such a way that not two of them are 
complementary.
These will be the $c_i 's$ and they have to be so that there
exists any way to put them in pairs, $\gamma = \{ (c_{\tau(1)},c_{\tau(2)}),
 \dots,(c_{\tau(n-1)},c_{\tau(n)})\}$ verifying
$$\sum_{(x,y)\in\gamma}d_{xy}+d_{v}=0.$$
Like in the precedent case, we assert that the border obtained is magic,
regardless the order in which $b_i's$ and $c_i 's$ have been placed in their
 respective rows or columns.
 
 The following theorem expresses the fact that this choice of the 
 numbers leads to a magic border.
\bigskip

 \begin {theorem}
 Let $A=\{v,w\}\cup\{b_{1},\dots,b_{n}\}\cup\{c_{1},\dots,c_{n}\}$ 
 and $B=A\cup A'$ a 
 border, built as we have just seen, then B is a magic border.
\end{theorem}

 \noindent{\bf Proof}
 
 ---Case $n$ even:

 Each row and column of the border must sum the magic constant
 of a square of order $n+2$, that is $S_{n+2}=\frac {(n+2)}{2}((n+2)^2+1)$;
 and we now that any pair of complementary numbers sum
 $(n+2)^2+1$, so in each row and column of the border we have the
 same sum as if we had $\frac {(n+2)}{2}$ pairs of complementary
 numbers. With the election done, we have
 $$\sum_{(x,y)\in\beta}d_{xy}=0$$so the total difference between these
 $\frac {(n+2)}{2}$ pairs and $\frac {(n+2)}{2}$ pairs of
 complementary numbers is zero, so for the rows of the border, the
 magic condition stands (we have seen it just for one
 row, but as the other is formed by the complementary numbers of
 our $b_i's$ and those of $v,w$ the condition stands for the other
 row).

 Now we have one column formed by $\{ c_1,\dots ,c_n \}\cup
 \{v,\bar w\}$ and we have $$\sum_{(x,y)\in\{ c_1, \dots ,c_n \}}
 d_{xy}=-d_{v\bar w}\Rightarrow \sum_{(x,y)\in\{ c_1, \dots ,c_n \}}
 d_{xy}+d_{v\bar w}=0$$ so columns also verify magic condition and
 hence, the border is magic.

---Case $n$ odd:

Now we must have in each row and column of the border the
 same sum as if we had $\frac {(n+1)}{2}$ pairs of complementary
 numbers, plus a half of the sum of one of these pairs, that is $\frac
 {(n+2)^2+1}{2}$. With the election done, we have
 $$v+\sum_{i=1}^nb_{i}+w=v+\frac {(n+1)}{2}((n+2)^2+1)+\sum_{(x,y)\in\beta}d_{xy}
 $$ $$=v+\frac {(n+1)}{2}((n+2)^2+1)-d_{v}=v+\frac {(n+1)}{2}((n+2)^2+1)-v+
 \frac{(n+2)^2+1}{2}$$ $$=\frac {(n+1)}{2}((n+2)^2+1)+\frac {((n+2)^2+1)}{2}$$
  so for the rows of the border, the
 magic condition stands.

 Now we have one column formed by$\{ c_1, \dots ,c_n \}\cup \{v,\bar w\}$
 and we have $$v+\sum_{i=1}^n c_{i}+\bar w=v+\frac {(n+1)}{2}((n+2)^2+1)+
 \sum_{(x,y)\in\gamma}d_{xy}$$ $$=v+\frac {(n+1)}{2}((n+2)^2+1)-d_{v}
 =v+\frac {(n+1)}{2}((n+2)^2+1)-v+
 \frac{(n+2)^2+1}{2}$$ $$=\frac {(n+1)}{2}((n+2)^2+1)+\frac {((n+2)^2+1)}{2},$$
  so columns also verify magic condition and, hence, the border is 
  magic. $\square$
\bigskip

Now we know that the choice we presented before lead us to a magic 
border. The next question is if that kind of choice is possible for 
any $n$. If so, we have a method for constructing magic squares of 
any order. The answer of this question is in the next theorem:

\begin {theorem}
 For any positive integer $n>2$, it is possible to make a choice of 
 the kind we have already seen , in order to build a magic border.
\end{theorem}
 \noindent {\bf Proof}
 
 The proof will consist on a constructive method. We will give an 
 actual example of choice verifying the needed properties. Some 
 generalizations of this method are possible.

  ---Case $n$ even:

  The way to prove this is to give (in a simple way) some election
  (among the various valid) that works for any $n$ even. We will
  proceed giving one for $n=4k,k\geq 1$ and other for
  $n=2(2k+1),k\geq 1$.

  $\cdot$ Case $n=4k,k\geq 1$

 We have the following arrangement of the numbers in two columns:

 \begin {picture}(180,180)(-80,20) \put(0,180){$1$}
\put(10,180){$\bullet$}
 \put(40,180){$\bullet$} \put(50,180){$(n+2)^2$}
 \put(0,165){$2$} \put(10,165){$\bullet$}
 \put(40,165){$\bullet$} \put(50,165){$(n+2)^2-1$}
 \put(0,150){$3$} \put(10,150){$\bullet$}
 \put(40,150){$\bullet$} \put(50,150){$(n+2)^2-2$}
 \put(0,135){$4$} \put(10,135){$\bullet$}
 \put(40,135){$\bullet$} \put(50,135){$(n+2)^2-3$}
 \put(0,120){$5$} \put(10,120){$\bullet$}
 \put(40,120){$\bullet$} \put(50,120){$(n+2)^2-4$}
 \put(0,105){$6$} \put(10,105){$\bullet$}
 \put(40,105){$\bullet$} \put(50,105){$(n+2)^2-5$}
 \put(0,90){$7$} \put(10,90){$\bullet$}
 \put(40,90){$\bullet$} \put(50,90){$(n+2)^2-6$}
 \put(0,75){$8$} \put(10,75){$\bullet$}
 \put(40,75){$\bullet$} \put(50,75){$(n+2)^2-7$}
 \put(0,60){$9$} \put(10,60){$\bullet$}
 \put(40,60){$\bullet$} \put(50,60){$(n+2)^2-8$}
 \put(0,45){\makebox(0,0)[br]{$10$}} \put(10,45){$\bullet$}
 \put(40,45){$\bullet$} \put(50,45){$(n+2)^2-9$}
 \put(12,30){$\vdots$}\put(42,30){$\vdots$}
\put(0,15){\makebox(0,15)[br]{$2n+2$}}\put(10,15){$\bullet$}
 \put(40,15){$\bullet$} \put(50,15){$n^2+2n+3$}
 \end{picture}

\bigskip
\noindent  where $(2n+2) \in\{10,18,26,\dots\}$.

 Now, we make our choice as follows: In a first part of the choice, we 
 take:
 $$v=(n+2)^2 -1,\qquad \qquad \bar w=5,$$
 $$b_1=1, \qquad b_2=3,\qquad b_3=(n+2)^2-3, \qquad b_4=(n+2)^2-4,$$
 $$c_1=6, \qquad c_2=(n+2)^2-7, \qquad c_3=9,\qquad c_4=(n+2)^2-9.$$
 Once this part of the choice is made (and we have finished if
 $n=4$), we have:
 $$d_{v\bar w}=(n+2)^2-1+5-(n^2+4n+5)=3,$$
 $$d_{vb_1}+d_{b_2b_3}+d_{b_4w}=-1-1+2=0,$$
 $$d_{c_1c_2}+d_{c_3c_4}=-3.$$
 And so, the magic property of the border stands in this first
 part.
 
 Now we have to place the other $2n-8$ numbers and its
 complementaries; $n-4$ of them will be $b_i's$ and the other $n-4$
 will be $c_i's$. Note that $n-4=4k-4=4(k-1)$, so we can divide
 these numbers into $k-1$ sets of $4$ numbers that we must place in a
 way in which the magic property keeps standing. An easy way of
 doing it is choosing number four by four in one of the following
 dispositions

  
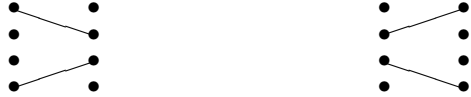
\begin {picture}(160,60)(-80,15)
\put(10,50){$\bullet$}
 \put(40,50){$\bullet$}
  \put(10,40){$\bullet$}
 \put(40,40){$\bullet$}
 \put(10,30){$\bullet$}
 \put(40,30){$\bullet$}
  \put(10,20){$\bullet$}
 \put(40,20){$\bullet$}
 \put(12,52){\line(3,-1){30}}
 \put(12,22){\line(3,1){30}}
 \put(150,50){$\bullet$}
 \put(180,50){$\bullet$}
  \put(150,40){$\bullet$}
 \put(180,40){$\bullet$}
 \put(150,30){$\bullet$}
 \put(180,30){$\bullet$}
  \put(150,20){$\bullet$}
 \put(180,20){$\bullet$}
\put(152,42){\line(3,1){30}}
 \put(152,32){\line(3,-1){30}} \end{picture}

\noindent or equivalently

 
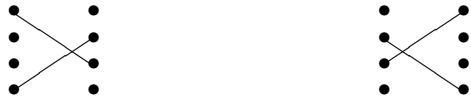
\begin {picture}(160,60)(-80,15)
\put(10,50){$\bullet$}
 \put(40,50){$\bullet$}
  \put(10,40){$\bullet$}
 \put(40,40){$\bullet$}
 \put(10,30){$\bullet$}
 \put(40,30){$\bullet$}
  \put(10,20){$\bullet$}
 \put(40,20){$\bullet$}
 \put(12,52){\line(3,-2){30}}
 \put(12,22){\line(3,2){30}}
 \put(150,50){$\bullet$}
 \put(180,50){$\bullet$}
  \put(150,40){$\bullet$}
 \put(180,40){$\bullet$}
 \put(150,30){$\bullet$}
 \put(180,30){$\bullet$}
  \put(150,20){$\bullet$}
 \put(180,20){$\bullet$}
\put(152,42){\line(3,-2){30}}
 \put(152,32){\line(3,2){30}} \end{picture}

 \noindent so, if we choose this way the remainding $b_i's$ and $c_i's$, we
 obtain a magic border (of course, once $b_i's$ have been
 chosen, we can put them in the border in any order, and the same
 stands for the $c_i's$, what makes the two dispositions above equivalent).

\begin{example}
 For example, we can form a magic border for $n=8$ with the
 numbers $1,2,\dots ,18$ and $83,84,\dots ,100$:

\bigskip
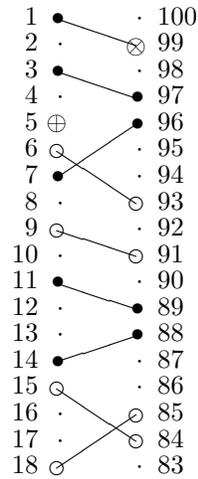
\begin {picture}(180,180)(-80,0) \put(0,170){$1$}
\put(10,170){$\bullet$}
 \put(42,170){$\cdot$} \put(50,170){$100$}
 \put(0,160){$2$} \put(12,160){$\cdot$}
 \put(38,159){$\otimes$}\put(50,160){$99$}
 \put(0,150){$3$} \put(10,150){$\bullet$}
 \put(42,150){$\cdot$} \put(50,150){$98$}
 \put(0,140){$4$} \put(12,140){$\cdot$}
 \put(40,140){$\bullet$} \put(50,140){$97$}
 \put(0,130){$5$} \put(8,130){$\oplus$}
 \put(40,130){$\bullet$} \put(50,130){$96$}
 \put(0,120){$6$} \put(12,122){\cir}
 \put(42,120){$\cdot$} \put(50,120){$95$}
 \put(0,110){$7$} \put(10,110){$\bullet$}
 \put(42,110){$\cdot$} \put(50,110){$94$}
 \put(0,100){$8$} \put(12,100){$\cdot$}
 \put(42,102){\cir} \put(50,100){$93$}
 \put(0,90){$9$} \put(12,92){\cir}
 \put(42,90){$\cdot$} \put(50,90){$92$}
 \put(0,80){\makebox(5,0)[br]{$10$}} \put(12,80){$\cdot$}
 \put(42,82){\cir} \put(50,80){$91$}
\put(0,70){\makebox(5,0)[br]{$11$}} \put(10,70){$\bullet$}
 \put(42,70){$\cdot$} \put(50,70){$90$}
 \put(0,60){\makebox(5,0)[br]{$12$}} \put(12,60){$\cdot$}
 \put(40,60){$\bullet$} \put(50,60){$89$}
 \put(0,50){\makebox(5,0)[br]{$13$}} \put(12,50){$\cdot$}
 \put(40,50){$\bullet$} \put(50,50){$88$}
 \put(0,40){\makebox(5,0)[br]{$14$}} \put(10,40){$\bullet$}
 \put(42,40){$\cdot$} \put(50,40){$87$}
 \put(0,30){\makebox(5,0)[br]{$15$}} \put(12,32){\cir}
 \put(42,30){$\cdot$} \put(50,30){$86$}
 \put(0,20){\makebox(5,0)[br]{$16$}} \put(12,20){$\cdot$}
 \put(42,22){\cir} \put(50,20){$85$}
 \put(0,10){\makebox(5,0)[br]{$17$}} \put(12,10){$\cdot$}
 \put(42,12){\cir} \put(50,10){$84$}
 \put(0,0){\makebox(5,0)[br]{$18$}} \put(12,2){\cir}
 \put(42,0){$\cdot$} \put(50,0){$83$}
\put(12,172){\line(3,-1){30}} \put(12,152){\line(3,-1){30}}
\put(12,122){\line(3,-2){30}}\put(12,112){\line(3,2){30}}
\put(12,92){\line(3,-1){30}}\put(12,72){\line(3,-1){30}}
\put(12,42){\line(3,1){30}}\put(12,32){\line(3,-2){30}}
\put(12,2){\line(3,2){30}} \put(100,110){where:}
\put(125,100){$\bullet\rightarrow b_i$} \put(128,93){\cir}\put
(132,91){ $\rightarrow c_i$} \put(125,80){$\otimes\rightarrow v$}
\put(125,70){$\oplus\rightarrow \bar w$}
 \end{picture}

\bigskip
 That produces the border
 $$\matrix{99&1&3&97&7&11&89&14&88&96 \cr
6&\quad&\quad&\quad&\quad&\quad&\quad&\quad&\quad&95\cr
93&\quad&\quad&\quad&\quad&\quad&\quad&\quad&\quad&8\cr
9&\quad&\quad&\quad&\quad&\quad&\quad&\quad&\quad&92\cr
91&\quad&\quad&\quad&\quad&\quad&\quad&\quad&\quad&10\cr
15&\quad&\quad&\quad&\quad&\quad&\quad&\quad&\quad&86\cr
85&\quad&\quad&\quad&\quad&\quad&\quad&\quad&\quad&16\cr
84&\quad&\quad&\quad&\quad&\quad&\quad&\quad&\quad&17\cr
18&\quad&\quad&\quad&\quad&\quad&\quad&\quad&\quad&83\cr
5&100&98&4&94&90&12&87&13&2}$$
\end{example}
\bigskip
 $\cdot$ Case $n=2(2k+1),k\geq 1$.

In this case, we have a similar arrangement of the numbers in two
columns so, making analogous considerations to those of the
precedent case, and noting that $2n+2\in\{14,22,\dots\}$, we can
make the following first part of the choice:
 $$v=1,\qquad \qquad \bar w=(n+2)^2-3,$$
 $$b_1=(n+2)^2-1, \qquad b_2=(n+2)^2-2,\qquad b_3=5,$$ $$ b_4=(n+2)^2-5,
 \qquad b_5=(n+2)^2-6,\qquad b_6=8,$$
 $$c_1=10, \qquad c_2=(n+2)^2-8, \qquad c_3=(n+2)^2-10,$$ $$ c_4=12,
 \qquad c_5=(n+2)^2-12,\qquad c_6=14,$$
 And for the second part, we proceed as before.

\begin{example}
 As an example, we can form a magic border for $n=10$ with the
 numbers $1,2,\dots ,22$ and $123,124,\dots ,144$:\bigskip

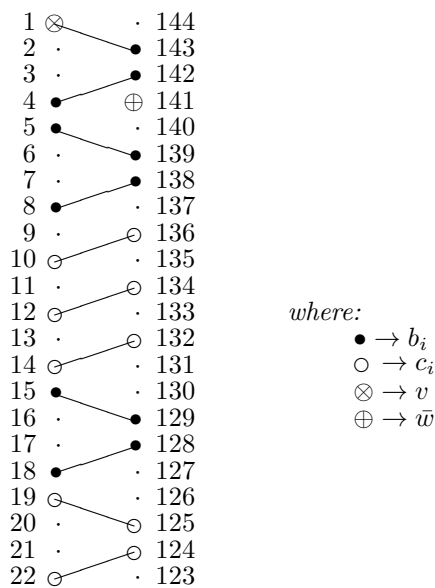
\begin {picture}(220,220)(-80,20) \put(0,220){$1$}
\put(8,220){$\otimes$}
 \put(42,220){$\cdot$} \put(50,220){$144$}
 \put(0,210){$2$} \put(12,210){$\cdot$}
 \put(40,210){$\bullet$}\put(50,210){$143$}
 \put(0,200){$3$} \put(12,200){$\cdot$}
 \put(40,200){$\bullet$} \put(50,200){$142$}
 \put(0,190){$4$} \put(10,190){$\bullet$}
 \put(38,190){$\oplus$} \put(50,190){$141$}
 \put(0,180){$5$} \put(10,180){$\bullet$}
 \put(42,180){$\cdot$} \put(50,180){$140$}
 \put(0,170){$6$} \put(12,170){$\cdot$}
 \put(40,170){$\bullet$} \put(50,170){$139$}
 \put(0,160){$7$} \put(12,160){$\cdot$}
 \put(40,160){$\bullet$} \put(50,160){$138$}
 \put(0,150){$8$} \put(10,150){$\bullet$}
 \put(42,150){$\cdot$} \put(50,150){$137$}
 \put(0,140){$9$} \put(12,140){$\cdot$}
 \put(42,142){\cir} \put(50,140){$136$}
 \put(0,130){\makebox(5,0)[br]{$10$}} \put(12,132){\cir}
 \put(42,130){$\cdot$} \put(50,130){$135$}
\put(0,120){\makebox(5,0)[br]{$11$}} \put(12,120){$\cdot$}
 \put(42,122){\cir} \put(50,120){$134$}
 \put(0,110){\makebox(5,0)[br]{$12$}} \put(12,112){\cir}
 \put(42,110){$\cdot$} \put(50,110){$133$}
 \put(0,100){\makebox(5,0)[br]{$13$}} \put(12,100){$\cdot$}
 \put(42,102){\cir} \put(50,100){$132$}
 \put(0,90){\makebox(5,0)[br]{$14$}} \put(12,92){\cir}
 \put(42,90){$\cdot$} \put(50,90){$131$}
 \put(0,80){\makebox(5,0)[br]{$15$}} \put(10,80){$\bullet$}
 \put(42,80){$\cdot$} \put(50,80){$130$}
 \put(0,70){\makebox(5,0)[br]{$16$}} \put(12,70){$\cdot$}
 \put(40,70){$\bullet$} \put(50,70){$129$}
 \put(0,60){\makebox(5,0)[br]{$17$}} \put(12,60){$\cdot$}
 \put(40,60){$\bullet$} \put(50,60){$128$}
 \put(0,50){\makebox(5,0)[br]{$18$}} \put(10,50){$\bullet$}
 \put(42,50){$\cdot$} \put(50,50){$127$}
 \put(0,40){\makebox(5,0)[br]{$19$}} \put(12,42){\cir}
 \put(42,40){$\cdot$} \put(50,40){$126$}
 \put(0,30){\makebox(5,0)[br]{$20$}} \put(12,30){$\cdot$}
 \put(42,32){\cir} \put(50,30){$125$}
 \put(0,20){\makebox(5,0)[br]{$21$}} \put(12,20){$\cdot$}
 \put(42,22){\cir} \put(50,20){$124$}
 \put(0,10){\makebox(5,0)[br]{$22$}} \put(12,12){\cir}
 \put(42,10){$\cdot$} \put(50,10){$123$}

\put(12,222){\line(3,-1){30}} \put(12,192){\line(3,1){30}}
\put(12,182){\line(3,-1){30}}\put(12,152){\line(3,1){30}}
\put(12,132){\line(3,1){30}}\put(12,112){\line(3,1){30}}
\put(12,92){\line(3,1){30}}\put(12,82){\line(3,-1){30}}
\put(12,52){\line(3,1){30}} \put(12,42){\line(3,-1){30}}
\put(12,12){\line(3,1){30}}
\put(100,110){where:}
\put(125,100){$\bullet\rightarrow b_i$} \put(128,93){\cir}\put
(132,91){ $\rightarrow c_i$} \put(125,80){$\otimes\rightarrow v$}
\put(125,70){$\oplus\rightarrow \bar w$}
 \end{picture}

\bigskip\bigskip
 That produces the border
 $$\matrix{1&143&142&5&139&138&8&15&129&128&18&4 \cr
136&\quad&\quad&\quad&\quad&\quad&\quad&\quad&\quad&\quad&\quad&9\cr
10&\quad&\quad&\quad&\quad&\quad&\quad&\quad&\quad&\quad&\quad&135\cr
134&\quad&\quad&\quad&\quad&\quad&\quad&\quad&\quad&\quad&\quad&11\cr
12&\quad&\quad&\quad&\quad&\quad&\quad&\quad&\quad&\quad&\quad&133\cr
132&\quad&\quad&\quad&\quad&\quad&\quad&\quad&\quad&\quad&\quad&13\cr
14&\quad&\quad&\quad&\quad&\quad&\quad&\quad&\quad&\quad&\quad&131\cr
19&\quad&\quad&\quad&\quad&\quad&\quad&\quad&\quad&\quad&\quad&126\cr
125&\quad&\quad&\quad&\quad&\quad&\quad&\quad&\quad&\quad&\quad&20\cr
124&\quad&\quad&\quad&\quad&\quad&\quad&\quad&\quad&\quad&\quad&21\cr
22&\quad&\quad&\quad&\quad&\quad&\quad&\quad&\quad&\quad&\quad&123\cr
141&2&3&140&6&7&137&130&16&17&127&144}$$
\end{example}
\bigskip

 ---Case $n$ odd:

 In this case we shall proceed in a somewhat different way. In this
 case, the method is general for all $n$ odd, but there are some
 details to observe. As in the even case, there are many ways to fill
 the border, among which we choose a simple one.
 We part from an arrangement of the $n+2$ pairs of complementary
 numbers as seen before. But now, we divide this columns in three
 parts:

 We call $A$ to the first part, containing the first $n-3$ pairs, that
 is $1 \rightarrow (n+2)^2$ to $n-3\rightarrow n^2+3n+8$.

 Part $B$ contains seven pairs, from $n-2 \rightarrow n^2+3n+7$ to
 $n+4 \rightarrow n^2+3n+1$.

 Finally, part $C$ contains the remaining $n-2$ pairs.

 Once we have these three sets, we select a number for the corner
 named $v$. We shall select   the number $(
 n+7)$ and put it in $v$ (we must note that this  $(n+7)$
  must be between $(n+5)$ and $(2n+2)$, and this happens if
 $n\geq 5$ so we must consider the case $n=3$ special. The reader 
 can construct a special method for $n=3$ easily from the general 
 method for $n$ odd). After
 this, we name $b_{i}'s$ or $c_{i}'s$ the numbers in the first column of part $C$ and in the
 second column of part $A$. Beginning with part $C$ we have two
 numbers above $(n+7)$, that is $(n+5)$ and $(n+6)$; we can name one of
 them $b_{i}$ for some $i$, and the other $c_{j}$ for some $j$,
 we match them with $(n+2)^2$ and $(n+2)^2-1$ respectively (labeled
 with the same letter of its pair). So we
 have that each of these pairs contributes with $(n+4)$ to $\sum d_{xy}'s$.
  The $(n-5)$ numbers under $v$ will be matched with
 the $n-5$ numbers under $(n+2)^2-1$, half of these
 pairs (note $(n-5)$ is even) being
 $b_{i}'s$ and the other half $c_{i}'s$ in any order. Each of these
 $n-5$ pairs contributes with $n+5$ to  $\sum d_{xy}$.

 Finally, numbers in set $B$ are matched following this scheme:

  
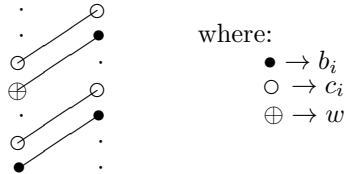
\begin {picture}(160,60)(-80,15)
  \put(12,60){$\cdot$}
 \put(42,62){\cir}
\put(12,50){$\cdot$}
 \put(40,50){$\bullet$}
  \put(12,42){\cir}
 \put(42,40){$\cdot$}
 \put(8,29){$\oplus$}
 \put(42,32){\cir}
  \put(12,20){$\cdot$}
 \put(40,20){$\bullet$}
 \put(12,12){\cir}
 \put(42,10){$\cdot$}
 \put(10,0){$\bullet$}
 \put(42,0){$\cdot$}
 \put(12,42){\line(3,2){30}}
 \put(12,32){\line(3,2){30}}
 \put(12,12){\line(3,2){30}}
 \put(12,2){\line(3,2){30}}
 \put(80,50){where:}
\put(105,40){$\bullet\rightarrow b_i$} \put(108,33){\cir}\put
(112,31){ $\rightarrow c_i$} \put(105,20){$\oplus\rightarrow w$}
  \end{picture}

 \bigskip\bigskip
 \noindent So, with this choice of the numbers, we must have

 $$\sum_{(x,y)\in\beta}d_{xy}+d_{v}=0$$
 $$\sum_{(x,y)\in\gamma}d_{xy}+d_{v}=0$$

 \noindent and we have $$d_{v}=(n+7)-(n^2+4n+5)=\frac{-(n^2+2n-9)}{2},$$

 $$\sum_{(x,y)\in\beta}d_{xy}=2(n+4)+\frac{(n-5)(n+5)}{2}+4=\frac{n^2+2n-9}{2},$$

 $$\sum_{(x,y)\in\gamma}d_{xy}=2(n+4)+\frac{(n-5)(n+5)}{2}+4=\frac{n^2+2n-9}{2}.$$

 Hence, our choice leads to a magic border for all $n$ odd.

\begin{example}
 As an example, we can form a magic border for $n=7$ with the
 numbers $1,2,\dots ,16$ and $66,67,\dots ,81$:

\bigskip
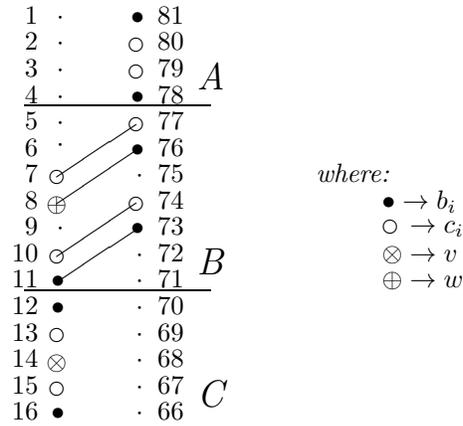
\begin {picture}(180,180)(-80,0) \put(0,170){$1$}
\put(12,170){$\cdot$}
 \put(40,170){$\bullet$} \put(50,170){$81$}
 \put(0,160){$2$} \put(12,160){$\cdot$}
 \put(42,162){\cir}\put(50,160){$80$}
 \put(0,150){$3$} \put(12,150){$\cdot$}
 \put (42,152){\cir} \put(50,150){$79$}
 \put(0,140){$4$} \put(12,140){$\cdot$}
 \put(40,140){$\bullet$} \put(50,140){$78$}
 \put(0,130){$5$} \put(12,130){$\cdot$}
 \put(42,132){\cir} \put(50,130){$77$}
 \put(0,120){$6$} \put(12,122){$\cdot$}
 \put(40,120){$\bullet$} \put(50,120){$76$}
 \put(0,110){$7$} \put(12,112){\cir}
 \put(42,110){$\cdot$} \put(50,110){$75$}
 \put(0,100){$8$} \put(8,99){$\oplus$}
 \put(42,102){\cir} \put(50,100){$74$}
 \put(0,90){$9$} \put(12,90){$\cdot$}
 \put(40,90){$\bullet$} \put(50,90){$73$}
 \put(0,80){\makebox(5,0)[br]{$10$}} \put(12,82){\cir}
 \put(42,80){$\cdot$} \put(50,80){$72$}
\put(0,70){\makebox(5,0)[br]{$11$}} \put(10,70){$\bullet$}
 \put(42,70){$\cdot$} \put(50,70){$71$}
 \put(0,60){\makebox(5,0)[br]{$12$}} \put(10,60){$\bullet$}
 \put(42,60){$\cdot$} \put(50,60){$70$}
 \put(0,50){\makebox(5,0)[br]{$13$}} \put(12,52){\cir}
 \put(42,50){$\cdot$} \put(50,50){$69$}
 \put(0,40){\makebox(5,0)[br]{$14$}} \put(8,39){$\otimes$}
 \put(42,40){$\cdot$} \put(50,40){$68$}
 \put(0,30){\makebox(5,0)[br]{$15$}} \put(12,32){\cir}
 \put(42,30){$\cdot$} \put(50,30){$67$}
 \put(0,20){\makebox(5,0)[br]{$16$}} \put(10,20){$\bullet$}
 \put(42,20){$\cdot$} \put(50,20){$66$}
\put(12,112){\line(3,2){30}} \put(12,102){\line(3,2){30}}
\put(12,82){\line(3,2){30}}\put(12,72){\line(3,2){30}}
\put(0,139){\line(1,0){70}}\put(0,69){\line(1,0){70}}
\put(65,145){\Large A}\put(65,75){\Large B}\put(65,25){\Large C}
\put(110,110){where:}
\put(135,100){$\bullet\rightarrow b_i$} \put(138,93){\cir}\put
(142,91){ $\rightarrow c_i$} \put(135,80){$\otimes\rightarrow v$}
\put(135,70){$\oplus\rightarrow w$}
 \end{picture}

\bigskip
 That produces the border
 $$\matrix{14&81&78&12&16&76&73&11&8 \cr
80&\quad&\quad&\quad&\quad&\quad&\quad&\quad&2\cr
79&\quad&\quad&\quad&\quad&\quad&\quad&\quad&3\cr
13&\quad&\quad&\quad&\quad&\quad&\quad&\quad&69\cr
15&\quad&\quad&\quad&\quad&\quad&\quad&\quad&67\cr
77&\quad&\quad&\quad&\quad&\quad&\quad&\quad&5\cr
7&\quad&\quad&\quad&\quad&\quad&\quad&\quad&75\cr
10&\quad&\quad&\quad&\quad&\quad&\quad&\quad&72\cr
74&1&4&70&66&6&9&71&68}$$
\end{example}
\bigskip

\begin {remark}
Let us insist in the fact that all the choices made in both even and odd cases were only particular ones, we have many possibilities that lead to other magic borders.
\end {remark}

\subsection {Generalization of the method.}

We can make a generalization of the method exposed in theorems 1 and 2 in three
ways (apart from those transformations of the square by its symetry
group):
\bigskip

\begin {remark}
Any permutation of the $b_{i}'s$ and $c_{i}'s$ in 
a magic border keeps the border magic.
\end {remark}

We have already said that $b_{i}'s$ and $c_{i}'s$ can be placed in
any order in their respective rows or columns when constructing the
magic border. Of course, when we change the disposition of  $b_{i}'s$
or $c_{i}'s$ we must change at the same time their respective
complementary numbers so that each pair of complementary numbers are
one opposite the other. As an example we can make a permutation of
both $b_{i}'s$ and $c_{i}'s$ in the magic border constructed in 
example 1:

 $$\matrix{99&88&14&89&11&7&97&3&1&96 \cr
93&\quad&\quad&\quad&\quad&\quad&\quad&\quad&\quad&8\cr
6&\quad&\quad&\quad&\quad&\quad&\quad&\quad&\quad&95\cr
91&\quad&\quad&\quad&\quad&\quad&\quad&\quad&\quad&10\cr
9&\quad&\quad&\quad&\quad&\quad&\quad&\quad&\quad&92\cr
85&\quad&\quad&\quad&\quad&\quad&\quad&\quad&\quad&16\cr
15&\quad&\quad&\quad&\quad&\quad&\quad&\quad&\quad&86\cr
18&\quad&\quad&\quad&\quad&\quad&\quad&\quad&\quad&83\cr
84&\quad&\quad&\quad&\quad&\quad&\quad&\quad&\quad&17\cr
5&13&87&12&90&94&4&98&100&2}$$
\bigskip

\noindent which is still magic. So, this method produces, if we consider
this generalitation, $2 \cdot n!$ magic borders with the same four numbers
in the corners.
\bigskip

\begin {remark}
 Changes in the construction of the inner
square do not affect the validity of the method.
\end {remark}

Changes in the construction of the inner square may have or not an
effect on our method. If we construct the inner square using
numbers $2n+3 \dots ,n^2+2n+2$ there is no effect on our method
whatever method we choose to construct this square. But if we use
different numbers in this construction, it will affect the
construction of the magic border(for example, we can use $\{1,\dots
\frac{n^2}{2}\}\cup\{(n+2)^2-\frac{n^2}{2}-1,\dots,(n+2)^2\}$). Anyway,
 we will have, for our border, $2n+1$ pairs of
complementary numbers. So, if we arrange them in two columns as we
saw before, we can apply the method in the same way, making the appropiate changes.
\bigskip

\begin {remark} There are many possibilities of changes in the 
scheme used in theorem 2 to verify the conditions of magic border 
that keep the magic of the border.
\end {remark}

The schemes used before in the proofs of the method for each case
are, as we said, only a valid choice among many correct ones. We
shall now have a closer look on this statement:

In any of the cases, we can say that the method consists of two
(more or less explicit) parts: in the first part we make the
choice of the corners, which leads to an {\it unbalance
situation}, and in the second part we try to {\it re-establish the
balance}, where the rest of the numbers in the border are
involved, being put by pairs of complementaries.

We have also seen in the different cases considered before two ways of solving this
{\it re-establishment} part: In the even cases, we used a minimum
number of pairs to re-establish the balance (first step of the
methods) and then, we pass to a second step that holds the {\it
inner balance} (second step of the methods). In the odd case we
have seen a way to re-establish the balance little by little
(which involved the pairs in $A$ and $C$) having a situation that
we can control, i.e. a {\it standard} step (scheme of the numbers
in $B$).

Any of these parts and steps can be changed, what leads us to
different methods (which are suitable to be generalised in the sense we have seen).
 We shall give now a pair of examples to illustrate these
affirmations.

\begin {example}
This example gives a magic border of order $5$ in which the
election of the corner is also changed, and pairs are not distributed
in three sets $A,B,C$.
\end {example}
\bigskip
\begin {picture}(120,120)(-40,0) \put(0,110){$1$}
\put(12,110){$\cdot$}
 \put(40,110){$\bullet$} \put(50,110){$49$}
 \put(0,100){$2$} \put(12,100){$\cdot$}
 \put(42,102){\cir}\put(50,100){$48$}
 \put(0,90){$3$} \put(12,90){$\cdot$}
 \put (42,92){\cir} \put(50,90){$47$}
 \put(0,80){$4$} \put(12,80){$\cdot$}
 \put(40,80){$\bullet$} \put(50,80){$46$}
 \put(0,70){$5$} \put(8,69){$\oplus$}
 \put(42,72){\cir} \put(50,70){$45$}
 \put(0,60){$6$} \put(12,62){\cir}
 \put(42,60){$\cdot$} \put(50,60){$44$}
 \put(0,50){$7$} \put(12,50){$\cdot$}
 \put(40,50){$\bullet$} \put(50,50){$43$}
 \put(0,40){$8$} \put(12,42){\cir}
 \put(42,40){$\cdot$} \put(50,40){$42$}
 \put(0,30){$9$} \put(10,30){$\bullet$}
 \put(42,30){$\cdot$} \put(50,30){$41$}
 \put(0,20){\makebox(5,0)[br]{$10$}} \put(12,22){\cir}
 \put(42,20){$\cdot$} \put(50,20){$40$}
\put(0,10){\makebox(5,0)[br]{$11$}} \put(8,8){$\otimes$}
 \put(42,10){$\cdot$} \put(50,10){$39$}
 \put(0,0){\makebox(5,0)[br]{$12$}} \put(10,0){$\bullet$}
 \put(42,00){$\cdot$} \put(50,0){$38$}
\put(12,72){\line(3,4){30}} \put(12,62){\line(3,4){30}}
\put(12,42){\line(3,5){30}}\put(12,32){\line(3,5){30}}
\put(12,22){\line(3,5){30}}\put(12,2){\line(3,5){30}}
\put (110,50) {$\matrix{11&9&12&49&46&43&5\cr
48&\quad&\quad&\quad&\quad&\quad&2\cr
47&\quad&\quad&\quad&\quad&\quad&3\cr
6&\quad&\quad&\quad&\quad&\quad&44\cr
8&\quad&\quad&\quad&\quad&\quad&42\cr
10&\quad&\quad&\quad&\quad&\quad&40\cr
45&41&38&1&4&7&39}$}
 \end{picture}

\begin {example}
This second example gives a magic border of order $8$ in which we have
changed the difference between the numbers in the corners, and the
distinction between the two parts of the method (in the even case).
\end {example}
 \bigskip
\begin {picture}(180,180)(-10,0) \put(0,170){$1$}
\put(12,170){$\cdot$}
 \put(40,170){$\bullet$} \put(50,170){$100$}
 \put(0,160){$2$} \put(10,160){$\bullet$}
 \put(42,160){$\cdot$}\put(50,160){$99$}
 \put(0,150){$3$} \put(12,150){$\cdot$}
 \put(40,150){$\bullet$} \put(50,150){$98$}
 \put(0,140){$4$} \put(10,140){$\bullet$}
 \put(42,140){$\cdot$} \put(50,140){$97$}
 \put(0,130){$5$} \put(12,132){$\cdot$}
 \put(42,132){\cir} \put(50,130){$96$}
 \put(0,120){$6$} \put(12,120){$\cdot$}
 \put (40,120){$\bullet$} \put(50,120){$95$}
 \put(0,110){$7$} \put(8,109){$\otimes$}
 \put(42,110){$\cdot$} \put(50,110){$94$}
 \put(0,100){$8$} \put(10,100){$\bullet$}
 \put (38,99){$\oplus$} \put(50,100){$93$}
 \put(0,90){$9$} \put(12,92){\cir}
 \put(42,90){$\cdot$} \put(50,90){$92$}
 \put(0,80){\makebox(5,0)[br]{$10$}} \put(12,80){$\cdot$}
 \put(40,80){$\bullet$} \put(50,80){$91$}
\put(0,70){\makebox(5,0)[br]{$11$}} \put(12,72){\cir}
 \put(42,70){$\cdot$} \put(50,70){$90$}
 \put(0,60){\makebox(5,0)[br]{$12$}} \put(10,60){$\bullet$}
 \put(42,60){$\cdot$} \put(50,60){$89$}
 \put(0,50){\makebox(5,0)[br]{$13$}} \put(12,50){$\cdot$}
 \put(40,50){$\bullet$} \put(50,50){$88$}
 \put(0,40){\makebox(5,0)[br]{$14$}} \put(12,40){$\cdot$}
 \put(42,42){\cir} \put(50,40){$87$}
 \put(0,30){\makebox(5,0)[br]{$15$}} \put(12,30){$\cdot$}
 \put(42,32){\cir} \put(50,30){$86$}
 \put(0,20){\makebox(5,0)[br]{$16$}} \put(12,22){\cir}
 \put(42,20){$\cdot$} \put(50,20){$85$}
 \put(0,10){\makebox(5,0)[br]{$17$}} \put(12,12){\cir}
 \put(42,10){$\cdot$} \put(50,10){$84$}
 \put(0,0){\makebox(5,0)[br]{$18$}} \put(12,0){$\cdot$}
 \put(42,2){\cir} \put(50,0){$83$}
\put(12,162){\line(3,1){30}} \put(12,142){\line(3,1){30}}
\put(12,112){\line(3,1){30}}\put(12,102){\line(3,-2){30}}
\put(12,92){\line(3,4){30}}\put(12,72){\line(3,-4){30}}
\put(12,62){\line(3,-1){30}}\put(12,22){\line(3,2){30}}
\put(12,12){\line(3,-1){30}}
\put(100,100){$\matrix{7&100&2&98&4&95&91&12&88&8 \cr
96&\quad&\quad&\quad&\quad&\quad&\quad&\quad&\quad&5\cr
9&\quad&\quad&\quad&\quad&\quad&\quad&\quad&\quad&92\cr
11&\quad&\quad&\quad&\quad&\quad&\quad&\quad&\quad&90\cr
87&\quad&\quad&\quad&\quad&\quad&\quad&\quad&\quad&14\cr
86&\quad&\quad&\quad&\quad&\quad&\quad&\quad&\quad&15\cr
16&\quad&\quad&\quad&\quad&\quad&\quad&\quad&\quad&85\cr
17&\quad&\quad&\quad&\quad&\quad&\quad&\quad&\quad&84\cr
83&\quad&\quad&\quad&\quad&\quad&\quad&\quad&\quad&18\cr
93&1&99&3&97&6&10&89&13&94}$}

 \end{picture}

\section {Even bordered magic squares: some results.}

In this section we give answer, with the help of the method we have just exposed, and for the even case, to some questions one may pose about magic borders.

We can see, from what we said in the {\it generalization of the method} that a magic border depends only on the number in the corners, and it is enough to consider two contiguous corners, say the two upper corners, which we denoted $v$ and $w$. So, a magic border can be represented as an element of $\Z^2\times \Z^n \times \Z^n$ being the first two numbers $\{v,w\}$, the second $n$ numbers the $b_i's$ and the last numbers the $c_i's$. We will denote by $\Omega_{v,w}^k$ the set of  magic borders of order $k$ having $v$ and $w$ as upper corners; $\Omega^k$ will stand for the set of magic borders of order $k$.

\begin {remark}\label{symmetries}
If we apply the symmetry group of the square to a given border, we obtain the following eight

$$\matrix{v&b_{1}&\dots&b_{n}&w \cr
c_{1}&\quad&\quad&\quad&\bar c_{1}\cr
\dots&\quad&\quad&\quad&\dots\cr
c_{n}&\quad&\quad&\quad&\bar c_{n}\cr
\bar w&\bar b_{1}&\dots&\bar b_{n}&\bar v}\qquad
\matrix{w&b_{n}&\dots&b_{1}&v \cr
\bar c_{1}&\quad&\quad&\quad&c_{1}\cr
\dots&\quad&\quad&\quad&\dots\cr
\bar c_{n}&\quad&\quad&\quad&c_{n}\cr
\bar v&\bar b_{1}&\dots&\bar b_{n}&\bar w}\qquad
\matrix{\bar w&\bar b_{1}&\dots&\bar b_{n}&\bar v \cr
c_{n}&\quad&\quad&\quad&\bar c_{n}\cr
\dots&\quad&\quad&\quad&\dots\cr
c_{1}&\quad&\quad&\quad&\bar c_{1}\cr
v&b_{1}&\dots&b_{n}&w}
$$
$$\matrix{\bar w&c_{n}&\dots&c_{1}&v \cr
\bar b_{1}&\quad&\quad&\quad&b_{1}\cr
\dots&\quad&\quad&\quad&\dots\cr
\bar b_{n}&\quad&\quad&\quad&b_{n}\cr
\bar v&\bar c_{n}&\dots&\bar c_{1}&w}\qquad
\matrix{\bar v&\bar b_{n}&\dots&\bar b_{1}&\bar w \cr
\bar c_{n}&\quad&\quad&\quad&c_{n}\cr
\dots&\quad&\quad&\quad&\dots\cr
\bar c_{1}&\quad&\quad&\quad&c_{1}\cr
w&b_{1}&\dots&b_{n}&v}\qquad
\matrix{w&\bar c_{1}&\dots&\bar c_{n}&\bar v \cr
b_{n}&\quad&\quad&\quad&\bar b_{n}\cr
\dots&\quad&\quad&\quad&\dots\cr
b_{1}&\quad&\quad&\quad&\bar b_{1}\cr
v&c_{1}&\dots&c_{n}&\bar w}
$$
$$\matrix{\bar v&\bar c_{n}&\dots&\bar c_{1}&w \cr
\bar b_{n}&\quad&\quad&\quad&b_{n}\cr
\dots&\quad&\quad&\quad&\dots\cr
\bar b_{1}&\quad&\quad&\quad&b_{1}\cr
\bar w&c_{n}&\dots&c_{1}&v}\qquad
\matrix{v&c_{1}&\dots&c_{n}&\bar w \cr
b_{1}&\quad&\quad&\quad&\bar b_{1}\cr
\dots&\quad&\quad&\quad&\dots\cr
b_{n}&\quad&\quad&\quad&\bar b_{n}\cr
w&\bar c_{1}&\dots&\bar c_{n}&\bar v}
$$
\bigskip
So if we know $\Omega ^k_{v,w}$ we can obtain, permuting its elements (considered in $\Z^{2n+2}$) all the elements in $\Omega ^k_{w,v}$, $\Omega ^k_{\bar w,\bar v}$, $\Omega ^k_{\bar w,v}$, $\Omega ^k_{\bar v,\bar w}$, $\Omega ^k_{w,\bar v}$, $\Omega ^k_{\bar v,w}$and $\Omega ^k_{v,\bar w}$.

\end {remark}

With this notation we can now set the following theorem:

\begin {theorem}\label {th3}
Let $n$ be a positive even number; let $v,w\in \{1,\dots, 2n+2\}$, then 
$$\Omega_{v,w}^n\neq \emptyset \iff \mbox {v and w have not the same parity}.$$
\end {theorem}

\noindent{\bf Proof:}

\noindent $\cdot$ We shall first prove the direct statement. To do this, we will see that if $v$ and $w$ have the same parity then $\Omega_{v,w}^n= \emptyset$:

We give here the proof for the case  $n=4k,\,k\in \N$; for the other case ($n=4k+2,\,k\in \N$) we can make analogous reasoning. 

If $n=4k$ then $S_n=2k((4k)^2+1)$ is even,  $S_{n+2}=(2k+1)((4k+2)^2+1)$ is odd, and $\sum_{i=1}^{2n+2}i=(2n+3)\frac{(2n+2)}{2}$ is odd.
If we arrange the numbers for the border in two columns (each number opposite its complementary) then we can take, with no loose of generality $v$ and $w$ from the first column.

Note that $\sum _i^nb_i=S_{4k+2}-v-w$ is odd, and $\sum _i^nc_i=S_{4k+2}-v-\bar w$ is even. 

Let $x$ be the sum of the $b_i's$ in the first column ( they are $\frac{n-2}{2}$ numbers) let $y$ be the sum of the $b_i's$ in the second column and let $x'$ and $y'$ be the sum of their respective complementaries. As $\sum _i^nb_i$ is odd, $x$ and $y$ cannot have the same parity;  and neither can $x$ and $x'$ nor $y$ and $y'$, because  $x'=\frac{n-2}{2}((n+2)^2+1)-x$, $y'=\frac{n+2}{2}((n+2)^2+1)-x$ being $\frac{n-2}{2}((n+2)^2+1)$ and $\frac{n+2}{2}((n+2)^2+1)$ odd numbers.

Now, let $t$ be the sum of the $c_i's$ in the first column ( they are $\frac{n}{2}$ numbers) let $s$ be the sum of the $c_i's$ in the second column and let $t'$ and $s'$ be the sum of their respective complementaries. As $\sum _i^nc_i$ is even, and so is $\frac{n}{2}((n+2)^2+1)$, we have $t,\,t',\,s$ and $s'$ have the same parity.

Finally, the sum of the numbers in the first column is $\sum_{i=1}^{2n+2}i$ which is odd, but this sum equals $v+w+x+y'+t+s'$ which, as we have just seen, is even, so if $v$ and $w$ have the same parity, then $\Omega_{v,w}^n= \emptyset$.

\noindent $\cdot$ Let us see now that if $v$ and $w$ have not the same parity, it is always possible to complete a magic border. 

First, note that, as we have already seen, it is enough to work the following cases: $$(1,2), \dots ,(1,2n+2)$$
$$ (2,3),\dots ,(2,2n+1)$$
$$ \dots \, \dots \, \dots$$
$$(2n,2n+1)$$
$$(2n+1,2n+2)$$

Note also that we can complete the borders of order $n+2$ by using those of order $n$, in the following way:

We use the numbers $1$ to $2n+2$ and their complementaries to complete the borders of order $4k$; now, for the order $4(k+1)$ we need $8$ new pairs of numbers. We shall match these numbers to have two pairs of $b_i's$ and two pairs of $c_i's$; let us do the matching so that $\sum d_{xy}=0$ for both the $b_i's$ and the $c_i's$; let us also avoid crossing lines in these matchings. We can make for instance, a matching like this:

 
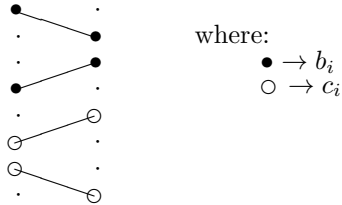
\begin {picture}(160,60)(-80,15)
  \put(10,60){$\bullet$}
 \put(42,60){$\cdot$}
\put(12,50){$\cdot$}
 \put(40,50){$\bullet$}
  \put(12,40){$\cdot$}
 \put(40,40){$\bullet$}
 \put(10,30){$\bullet$}
 \put(42,30){$\cdot$}
  \put(12,20){$\cdot$}
 \put(42,22){\cir}
 \put(12,12){\cir}
 \put(42,10){$\cdot$}
 \put(12,2){\cir}
 \put(42,0){$\cdot$}
 \put(12,-10){$\cdot$}
 \put(42,-8){\cir}
 \put(12,62){\line(3,-1){30}}
 \put(12,32){\line(3,1){30}}
 \put(12,12){\line(3,1){30}}
 \put(12,2){\line(3,-1){30}}
 \put(80,50){where:}
\put(105,40){$\bullet\rightarrow b_i$} \put(108,33){\cir}\put
(112,31){ $\rightarrow c_i$} 
  \end{picture}

 \bigskip
 
 \bigskip
 
 \bigskip
 Now, we can put this scheme at the bottom of the scheme of a border of order $4k$ with $v$ and $w$ as upper corners and we have a border of order $4(k+1)$ with the same upper corners. If we put the first two pairs at the top of the scheme of a magic border of order $n$ with $v$ and $w$ as upper corners, and the other six pairs at the bottom, we obtain a scheme for a magic border of order $n+2$ with $v+2$ and $w+2$ as upper corners. With the same procedure, starting from an element of $\Omega^{4k}_{v,w}$ we have elements of $\Omega^{4(k+1)}_{v,w}$, $\Omega^{4(k+1)}_{v+2,w+2}$, $\Omega^{4(k+1)}_{v+4,w+4}$, $\Omega^{4(k+1)}_{v+6,w+6}$ and $\Omega^{4(k+1)}_{v+8,w+8}$. 
 
 So, starting from the magic borders of order $4k$ we can obtain  elements of $\Omega^{4(k+1)}_{v,w}$ for all suitable pairs $(v,w)$, except for the following:
 $$(1,2m-4),(1,2m-2),(1,2m),(1,2m+2)$$
 $$(2,2m-5),(2,2m-3),(2,2m-1),(2,2m+1)$$
 $$(3,2m-2),(3,2m),(3,2m+2)$$
 $$(4,2m-3),(4,2m-1),(4,2m+1)$$
 $$(5,2m),(5,2m+2)$$
 $$(6,2m-1),(6,2m+1)$$
 $$(7,2m+2)$$
 $$(8,2m+1)$$
 where m=4(k+1).
 
  Thus, what we need is to verify if $\Omega^4_{v,w}\neq \emptyset$ for all suitable pairs $(v,w)$ and $\Omega^{4(k+1)}_{v,w}\neq \emptyset$ for the 20 pairs $(v,w)$ that cannot be solved starting from elements of $\Omega^{4(k)}$.
 
 This can be easily done using the method exposed in the first section, and this completes the proof. As it would be too long to draw here all the schemes,we give only the solutions, in two tables. In the first one, we write an element of $\Omega_{v,w}^4$ for every one of the 25 pairs needed, and in the second, an element of $\Omega_{v,w}^m$ for every one of the pairs needed when $m=4(k+1)$, being in this second table
 $$B= \bigcup_{i=1}^{\frac {m-8}{4}} \{11+8 i,16+8i,m^2+2m-2+8i,m^2+2m+1+8i\}$$
 $$C= \bigcup_{i=1}^{\frac {m-8}{4}} \{13+8 i,14+8i,m^2+2m+3+8i,m^2+2m+4+8i\}$$
 
 This completes the proof of our theorem.
 
 \bigskip
 \bigskip
 
 \begin {table}
 \begin {center}
  \begin {tabular}{|c c|c|c|}

 \hline 
 v&w&$b_i's$&$c_i's$\\
 \hline \hline
 1&2&34,33,32,9&6,30,29,10\\ \hline
 1&4&35,32,31,8&34,7,9,27\\ \hline
 1&6&35,34,5,30&33,8,28,10\\ \hline
 1&8&35,3,33,31&32,30,9,10\\ \hline
 1&10&35,4,31,30&34,32,8,9\\ \hline
 3&2&36,32,30,8&4,31,28,10\\ \hline
 3&4&36,35,5,28&6,30,29,10\\ \hline
 3&6&36,35,4,27&32,7,29,9\\ \hline
 3&8&36,5,31,28&35,4,30,10\\ \hline
 3&10&1,35,32,30&33,31,8,9\\ \hline
 5&2&36,34,4,30&6,29,9,27\\ \hline
 5&4&36,2,34,30&6,29,28,10\\ \hline
 5&6&36,2,34,28&33,7,8,27\\ \hline
 5&8&1,34,33,30&35,6,9,27\\ \hline
 5&10&1,34,33,28&35,6,30,8\\ \hline
 7&2&36,3,32,31&4,29,9,27\\ \hline
 7&4&36,2,34,28&5,31,8,27\\ \hline
 7&6&1,35,34,28&33,5,8,27\\ \hline
 7&8&36,5,28,27&2,34,33,6\\ \hline
 7&10&1,33,32,28&35,3,31,8\\ \hline
 9&2&36,3,32,29&4,6,30,27\\ \hline
 9&4&36,2,31,29&3,32,7,27\\ \hline
 9&6&1,35,33,27&3,32,7,29\\ \hline
 9&8&1,35,31,27&34,4,5,30\\ \hline
 9&10&1,32,30,29&2,34,33,6\\ 
 
 \hline
 \end {tabular}
 \end {center}
 \caption{$\Omega^4_{v,w}$}
 \end {table}

 \setlength{\oddsidemargin}{18pt}
 \begin {table}

 \begin {tabular}{|c c|p{6.5cm}| p{6.5cm}|}\hline 
 v&w&$b_i's=B\quad\cup$&$c_i's=C\quad\cup$\\
 \hline \hline
 1&2m+2&$\{10,14,16,m^2+2m+4,(m+2)^2-15, $&$ \{7,8,9,11,(m+2)^2-5,(m+2)^2-4,$
 \\
 & & $(m+2)^2-12,(m+2)^2-11,(m+2)^2-1\}$&$(m+2)^2-3,(m+2)^2-2\}$\\ \hline
  
 3&2m+2&$\{7,13,14,m^2+2m+4,(m+2)^2-15, $&$ \{4,9,10,15,(m+2)^2-10,(m+2)^2-5,$
 \\
  & & $(m+2)^2-11,(m+2)^2-7,(m+2)^2-1\}$& $(m+2)^2-4,(m+2)^2\}$\\ \hline
  
5&2m+2&$\{2,13,2m-1,m^2+2m+4,m^2+2m+5, $&$ \{4,10,11,14,(m+2)^2-11,(m+2)^2-6,$
 \\
  & & $(m+2)^2-8,(m+2)^2-7,(m+2)^2-2\}$& $(m+2)^2-5,(m+2)^2\}$\\ \hline
  
7&2m+2&$\{3,5,15,m^2+2m+4,(m+2)^2-10, $&$ \{2,12,13,16,(m+2)^2-13,(m+2)^2-8,$
 \\
  & & $(m+2)^2-9,(m+2)^2-5,(m+2)^2-3\}$& $(m+2)^2-7,(m+2)^2\}$\\ \hline
  
2&2m+1&$\{7,13,16,m^2+2m+3,(m+2)^2-13, $&$ \{4,9,10,15,(m+2)^2-10,(m+2)^2-5,$
 \\
  & & $(m+2)^2-11,(m+2)^2-7,(m+2)^2-2\}$& $(m+2)^2-4,(m+2)^2\}$\\ \hline
  
4&2m+1&$\{3,13,2m-1,m^2+2m+3,m^2+2m+5, $&$ \{5,10,11,14,(m+2)^2-11,(m+2)^2-6,$
 \\
  & & $(m+2)^2-8,(m+2)^2-7,(m+2)^2\}$& $(m+2)^2-5,(m+2)^2-1\}$\\ \hline
  
6&2m+1&$\{3,9,16,m^2+2m+3,(m+2)^2-13, $&$ \{2,11,12,15,(m+2)^2-12,(m+2)^2-7,$
 \\
  & & $(m+2)^2-9,(m+2)^2-4,(m+2)^2-3\}$& $(m+2)^2-6,(m+2)^2\}$\\ \hline
  
8&2m+1&$\{4,14,15,m^2+2m+3,(m+2)^2-15, $&$ \{2,5,12,13,(m+2)^2-9,(m+2)^2-8,$
 \\
  & & $(m+2)^2-10,(m+2)^2-6,(m+2)^2-5\}$& $(m+2)^2-2,(m+2)^2\}$\\ \hline
  
1&2m&$\{8,9,2m-1,m^2+2m+4,m^2+2m+8, $&$ \{5,10,11,2m+2,m^2+2m+7,$
 \\
  & & $(m+2)^2-11,(m+2)^2-3,(m+2)^2-2\}$& $(m+2)^2-6,(m+2)^2-5,(m+2)^2-1\}$\\ \hline
  
3&2m&$\{7,12,2m+2,m^2+2m+4,m^2+2m+6, $&$ \{4,9,10,13,(m+2)^2-10,(m+2)^2-5,$
 \\
  & & $(m+2)^2-13,(m+2)^2-7,(m+2)^2-1\}$& $(m+2)^2-4,(m+2)^2\}$\\ \hline
  
5&2m&$\{2,14,2m+1,m^2+2m+3,m^2+2m+6, $&$ \{3,10,11,13,(m+2)^2-11,(m+2)^2-6,$
 \\
  & & $(m+2)^2-12,(m+2)^2-11,(m+2)^2-1\}$& $(m+2)^2-5,(m+2)^2\}$\\ \hline
  
2&2m-1&$\{8,13,2m+2,m^2+2m+4,m^2+2m+5, $&$ \{3,9,10,14,(m+2)^2-10,(m+2)^2-5,$
 \\
  & & $(m+2)^2-11,(m+2)^2-6,(m+2)^2-3\}$& $(m+2)^2-4,(m+2)^2\}$\\ \hline
  
4&2m-1&$\{8,12,2m+2,m^2+2m+4,m^2+2m+5, $&$ \{2,9,10,13,(m+2)^2-10,(m+2)^2-5,$
 \\
  & & $(m+2)^2-13,(m+2)^2-6,(m+2)^2-2\}$& $(m+2)^2-4,(m+2)^2\}$\\ \hline
  
6&2m-1&$\{3,11,14,m^2+2m+4,(m+2)^2-12, $&$ \{4,9,10,2m+2,(m+2)^2+2m+5,(m+2)^2-7,$
 \\
  & & $(m+2)^2-11,(m+2)^2-4,(m+2)^2-1\}$& $(m+2)^2-6,(m+2)^2\}$\\ \hline
  
1&2m-2&$\{5,10,2m,m^2+2m+3,m^2+2m+8, $&$ \{6,11,12,2m+1,m^2+2m+6,(m+2)^2-7,$
 \\
  & & $(m+2)^2-8,(m+2)^2-3,(m+2)^2-1\}$& $(m+2)^2-6,(m+2)^2-2\}$\\ \hline
  
3&2m-2&$\{7,2m-1,2m+2,m^2+2m+4,m^2+2m+5, $&$ \{2,9,10,2m-3,m^2+2m+10,(m+2)^2-5,$
 \\
  & & $m^2+2m+9,(m+2)^2-7,(m+2)^2-3\}$& $(m+2)^2-4,(m+2)^2\}$\\ \hline
  
2&2m-3&$\{7,11,2m+2,m^2+2m+5,m^2+2m+7, $&$ \{4,9,10,2m+1,m^2+2m+6,(m+2)^2-5,$
 \\
  & & $(m+2)^2-11,(m+2)^2-7,(m+2)^2\}$& $(m+2)^2-4,(m+2)^2-2\}$\\ \hline
  
4&2m-3&$\{11,2m-1,2m+1,m^2+2m+3,m^2+2m+5, $&$ \{3,7,8,12,(m+2)^2-8,(m+2)^2-5,$
 \\
  & & $m^2+2m+7,(m+2)^2-11,(m+2)^2-1\}$& $(m+2)^2-4,(m+2)^2\}$\\ \hline
  
1&2m-4&$\{2m-3,2m,2m+1,m^2+2m+3,m^2+2m+6, $&$ \{6,7,8,2m-5,m^2+2m+12,(m+2)^2-4,$
 \\
  & & $m^2+2m+7,m^2+2m+11,(m+2)^2-1\}$& $(m+2)^2-3,(m+2)^2-2\}$\\ \hline
  
2&2m-5&$\{8,2m,2m+2,m^2+2m+4,m^2+2m+7, $&$ \{5,6,9,2m-1,m^2+2m+9,(m+2)^2-6,$
 \\
  & & $m^2+2m+8,(m+2)^2-11,(m+2)^2-1\}$& $(m+2)^2-3,(m+2)^2-2\}$\\ 
 \hline
 \end {tabular}

 \caption {$\Omega^m_{v,w};\, m=4(k+1)$ for the missing pairs.}
 \end{table}

 \section {Conclusions.}
 A method for the construction of magic bordered squares was presented, which consists of a simplification of the problem of construcing magic borders. The method allows the construction of such borders in an easy way by arranging numbers in a two-column diagram. The method was used to give a complete description of magic borders of even order. 
 
 In addition, the method can be used (among other utilities) to classify the existing methods for constructing magic squares (which are known since the 10th century) and to give a complete enumeration of magic bordered squares for some orders. This will be exposed in some other paper.

\bibliographystyle{amsplain}

 \end {document}